\documentclass[11pt,a4paper]{article}

\usepackage[T1]{fontenc}
\usepackage[utf8]{inputenc}
\usepackage{authblk}

\usepackage{amssymb}
\usepackage{amsmath}
\usepackage{pgfplots}
\usepackage{subfigure}
\usepackage{caption}
\usepackage{epsfig}
\usepackage{url}

\def\be{\begin{equation}}
\def\ee{\end{equation}}

\newtheorem{theorem}{Theorem}
 \newenvironment{keywords}{ \noindent {\small \textbf{Key Words}}.}{ }

\begin{document}


\newcommand{\nms}{\normalsize}

\title{A deterministic global optimization using smooth diagonal auxiliary functions
\thanks{This work was supported by the INdAM--GNCS 2014 Research
Project of the Italian National Group for Scientific Computation of
the National Institute for Advanced Mathematics ``F.\,Severi''.}}

\author[1,2]{ Ya.D. Sergeyev\thanks{ email: yaro@dimes.unical.it }\thanks{Corresponding author. }}
\author[1,2]{Dmitri E. Kvasov\thanks{email: kvadim@dimes.unical.it }
}
 \nms
 \affil[1]{\nms DIMES, Universit\`a della  Calabria, Italy}
\affil[2]{\nms Department of Software and Supercomputing
Technologies\\ Lobachevsky State University of Nizhni Novgorod,
Russia}

\date{}

\maketitle

\begin{abstract}
In many practical decision-making problems it happens that
functions involved in optimization process  are black-box with
unknown analytical representations and hard to evaluate. In this
paper, a global optimization problem is considered where both the
goal function~$f(x)$ and its gradient $f'(x)$ are black-box
functions. It is supposed that $f'(x)$ satisfies the Lipschitz
condition over the search hyperinterval with an unknown Lipschitz
constant~$K$. A new deterministic `Divide-the-Best' algorithm
based on efficient diagonal partitions and smooth auxiliary
functions is proposed in its basic version, its convergence
conditions are studied and numerical experiments executed on eight
hundred test functions are presented.
\end{abstract}



\section{Introduction}

In many important applied problems, some decisions should be made by
finding the global optimum of a multiextremal objective function
subject to a set of constrains (see, e.\,g.,
\cite{Floudas&Pardalos(2009), Grossmann(1996), Horst&Pardalos(1995),
Mockus(2000), Pardalos&Romeijn(2002),Paulavicius&Zilinskas(2014a),
Pinter(1996), Sergeyev&Kvasov(2008), Sergeyev&Kvasov(2011),
Strongin(1978), Strongin&Sergeyev(2000),
Zhigljavsky&Zilinskas(2008)} and the references given therein).
Frequently, especially in engineering applications, the functions
involved in optimization process are black-box with unknown
analytical representations and hard to evaluate. Such
computationally challenging decision-making problems often cannot be
solved by traditional optimization techniques based on strong
suppositions (such as, e.\,g., convexity) about the problem.

Because of significant computational costs involved, normally a
limited number of functions evaluations are available for a
decision-maker (physicist, biologist, economist, engineer, etc.)
when he/she optimizes this kind of functions. Hence, the main goal
is to construct fast global optimization algorithms that generate
acceptable solutions with a relatively small number of functions
evaluations.

Usually, the methods applied in this context are subdivided in
stochastic and deterministic. Stochastic approaches (see, e.\,g.,
\cite{Calvin&Zilinskas(2000), Floudas&Pardalos(2009),
Horst&Pardalos(1995), Mockus(2000), Pardalos&Romeijn(2002),
Zhigljavsky(1991), Zhigljavsky&Zilinskas(2008)}) can often work in
a simpler manner than the deterministic algorithms. They can be
also suitable for the problems where the functions evaluations are
corrupted by noise. However, solutions found by some stochastic
algorithms (e.\,g., by popular heuristic methods like evolutionary
algorithms, simulated annealing, etc.; see, e.\,g.,
\cite{Floudas&Pardalos(2009), Gandomi:et:al.(2013),
Pardalos&Romeijn(2002), Schneider&Kirkpatrick(2006), Yang(2010)})
can be only local solutions that can be located far away from the
global one. This fact can preclude these methods from their usage
in practice if an accurate estimate of the global solution is
required. Therefore, our attention in the paper is focused on
deterministic approaches.

Deterministic global optimization is an important applied field
(see, e.\,g., \cite{Floudas&Pardalos(2009), Horst&Pardalos(1995),
Pinter(1996), Sergeyev&Kvasov(2008), Strongin&Sergeyev(2000),
Zilinskas&Zilinskas(2010)}). As a rule, deterministic methods
exhibit (under certain conditions) rigorous global convergence
properties and allow one to obtain guaranteed estimations of
global solutions. However, deterministic models can still require
too many function evaluations to achieve adequately good solutions
for the considered global optimization problems.

One of the simplest techniques in this framework is represented by
the so-called derivative-free (or direct) approach (see, e.\,g.,
\cite{Conn:et:al.(2009), Custodio:et:al.(2011),
Kolda:et:al.(2003), Liuzzi:et:al.(2010), Rios&Sahinidis(2013)}),
often used for solving important applied problems (see, e.\,g.,
pattern search methods~\cite{Kolda:et:al.(2003)}, the DIRECT
method \cite{Jones:et:al.(1993)}, the response surface, or
surrogate model methods~\cite{Jones:et:al.(1998)}, etc.).
Unfortunately (see, e.\,g.,~\cite{DiSerafino:et:al.(2011),
Kvasov&Sergeyev(2012a), Paulavicius&Zilinskas(2014b),
Sergeyev&Kvasov(2006)}), these methods either are designed to
detect only stationary points or can demand too high computational
effort for their work. As observed, e.\,g., in
\cite{Horst&Pardalos(1995), Stephens&Baritompa(1998)}, if no
particular assumptions are made on the objective function, any
finite number of function evaluations cannot guarantee getting
close to the global minimum, since the function may have very deep
and narrow peaks.

The Lipschitz-continuity assumption is one of the natural
suppositions on the objective function (in fact, in technical
systems the energy of change is always limited). The problem
involving Lipschitz functions is said to be the Lipschitz global
optimization problem (see, e.\,g.,~\cite{Evtushenko(1985),
Horst&Pardalos(1995), Kvasov&Sergeyev(2012b), Pinter(1996),
Sergeyev&Kvasov(2008), Sergeyev&Kvasov(2011),
Strongin&Sergeyev(2000), Zhigljavsky&Zilinskas(2008)}).

The global optimization problem with a differentiable objective
function having the Lipschitz gradient (with an unknown Lipschitz
constant) is an important class of Lipschitz global optimization
problems. Formally, this class of problems can be stated as
follows:

 \be \label{LGOP_f}
   f ^{*}=f(x^{*})=\min_{x \in D} \; f(x),
 \ee
 \be \label{LGOP_K}
   \| f'(x') - f'(x'') \| \le K \| x'- x'' \|,
   \hspace*{3mm}   x', x''\in D, \hspace*{3mm} 0 <K < \infty,
 \ee
where
 \be \label{LGOP_D}
   D = [a,b] = \{ x \in R^N: a(j) \leq x(j) \leq b(j) \}.
 \ee

It is assumed here that the objective function $f(x)$ can be
black-box, multiextremal, its gradient $f'(x)=\left(\frac{\partial
f(x)}{\partial x(1)}, \frac{\partial f(x)}{\partial x(2)}, \ldots,
\frac{\partial f(x)}{\partial x(N)}\right)^{T}$ (which can be
itself an expensive multiextremal black-box vector-function) can
be calculated during the search, and $f'(x)$ is
Lipschitz-continuous with some unknown constant~$K$, $0<K<\infty$,
over $D$. Problem (\ref{LGOP_f})--(\ref{LGOP_D}) is frequently met
in engineering (see, e.\,g., \cite{Pinter(1996),
Sergeyev&Kvasov(2008), Strongin&Sergeyev(2000)}), for instance, in
electrical engineering design (see, e.\,g.,
\cite{Sergeyev:et:al.(1999), Strongin&Sergeyev(2000)}).

There are known several methods for solving this problem that can
be distinguished with respect to the way the Lipschitz
constant~$K$ from (\ref{LGOP_D}) is estimated in their
computational schemes. There exist algorithms using an a priori
given estimate of $K$ (see, e.\,g., \cite{Breiman&Cutler(1993),
Sergeyev(1998a)}), its adaptive estimates (see, e.\,g.,
\cite{Gergel(1997), Gergel&Sergeyev(1999), Sergeyev(1998a)}), and
adaptive estimates of local Lipschitz constants (see, e.\,g.,
\cite{Sergeyev(1998a), Sergeyev&Kvasov(2008)}). Recently, methods
working with multiple estimates of $K$ chosen from a set of
possible values have been also proposed
(see~\cite{Kvasov&Sergeyev(2009), Kvasov&Sergeyev(2012b)}).

This paper is devoted to developing a new global optimization
algorithm for solving problem (\ref{LGOP_f})--(\ref{LGOP_D}). The
new method adaptively estimates the unknown Lipschitz constant $K$
from (\ref{LGOP_K}) during the search and is based on efficient
diagonal partitions of the search domain proposed recently.
Theoretical background for the introduced method is drawn in
Section 2. The algorithm is presented and analyzed theoretically
in Section~3. Section~4 presents results of numerical experiments
executed on several hundreds of multiextremal test functions.
Finally, Section~5 concludes the paper.

\section{Theoretical background}

In this Section, some important theoretical results, necessary for
introducing the new algorithm, are briefly described.

\subsection{`Divide-the-Best' algorithms}

As known, many global optimization methods (of both stochastic and
deterministic types) have a similar structure. Therefore, several
approaches to the development of a general framework for
describing global optimization algorithms and providing their
convergence conditions in a unified manner have been proposed
(see, e.\,g., \cite{Floudas&Pardalos(2009),
Grishagin:et:al.(1997), Horst&Pardalos(1995), Pinter(1996)}). The
`\emph{Divide-the-Best}' approach (see \cite{Sergeyev(1998b),
Sergeyev&Kvasov(2008)}) is one of such unifying schemes. It
generalizes both the schemes of adaptive
partition~\cite{Pinter(1996)} and
characteristic~\cite{Grishagin:et:al.(1997),
Sergeyev&Kvasov(2008), Strongin&Sergeyev(2000)} algorithms, which
are widely used for constructing and studying global optimization
methods.

\begin{figure}[ht]
\begin{center}
\includegraphics[width=100mm,keepaspectratio]{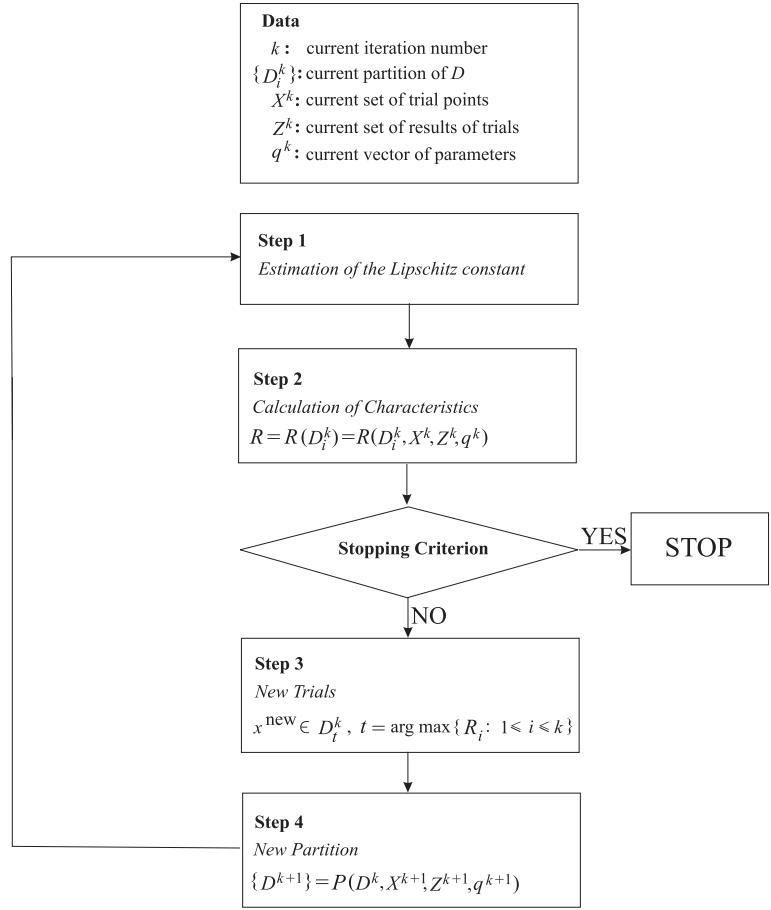}
\caption{Flow chart of `Divide-the-Best' algorithms.}
\label{fig:DBA}
\end{center}
\end{figure}

In a `Divide-the-Best' algorithm (its generic iteration is
represented by the flow chart in Figure~\ref{fig:DBA}), given a
vector $q$ of the method parameters, an adaptive partition of the
search domain $D$ from~(\ref{LGOP_D}) into a finite set
$\{D^k_i\}$ of robust subsets $D^k_i$ is considered at every
iteration $k$. In Step 1, the Lipschitz constant ($K$
from~(\ref{LGOP_K}), for the objective function gradient $f'(x)$;
or $L$ in the case of the Lipschitz objective function $f(x)$) is
estimated in some way. Basing on the previously obtained
information $X^k$, $Z^k$ about the objective function, the `merit'
(called \emph{characteristic}) $R_i$ of each subset (see Step~2 in
Figure~\ref{fig:DBA}) is estimated for performing a further, more
detailed, investigation (see Steps~3 and 4 in
Figure~\ref{fig:DBA}). The best (in a predefined sense)
characteristic achieved over a hyperinterval $D^k_t$ corresponds
to a higher possibility to determine the global minimum point
within $D^k_t$ (see Step 3). The hyperinterval $D^k_t$ is then
partitioned at the next iteration of the algorithm. More than one
`promising' hyperintervals can be subdivided at each iteration.

Various techniques (for instance, in the context of the
\emph{geometric} approach, see, e.\,g.,
\cite{Horst&Pardalos(1995), Kvasov&Sergeyev(2012b),
Sergeyev&Kvasov(2011), Strongin&Sergeyev(2000)}) for selection of
hyperintervals to be partitioned (see Step~3 in
Figure~\ref{fig:DBA}) and for executing this partitioning (by
means of an operator $P$, see Step~4 in Figure~\ref{fig:DBA} and
the following subsection 2.2) can be used within this scheme.

As the stopping criteria, one can check, e.\,g., the volume of a
hyperinterval with the best characteristic or depletion of
computational resources such as the maximal number of
\emph{trials} (i.\,e., evaluations of the objective function and
its gradient).

Theoretical analysis of the `Divide-the-Best' algorithms with
different types of characteristics and partition operators is
performed in~\cite{Sergeyev(1998b), Sergeyev&Kvasov(2008)}. A
particular attention is given there to situations (important in
applied problems) when conditions of global (local) convergence
are satisfied not in the whole domain~$D$, but only in its small
subregion(s). This can happen, e.\,g., in many Lipschitz global
optimization algorithms that either underestimate the Lipschitz
constant during their work or use local information in various
subregions of $D$ (see, e.\,g., \cite{Kvasov:et:al.(2003),
Sergeyev(1998b), Sergeyev&Kvasov(2008), Strongin&Sergeyev(2000)}).
It should be also noticed that the `Divide-the-Best' scheme can be
successfully applied to develop parallel multidimensional global
optimization methods (see, e.\,g., \cite{Grishagin:et:al.(1997),
Strongin&Sergeyev(2000)}).

\subsection{Efficient partition strategy}

\begin{figure}[!t]
\centering{\includegraphics[width=0.7\linewidth,keepaspectratio]{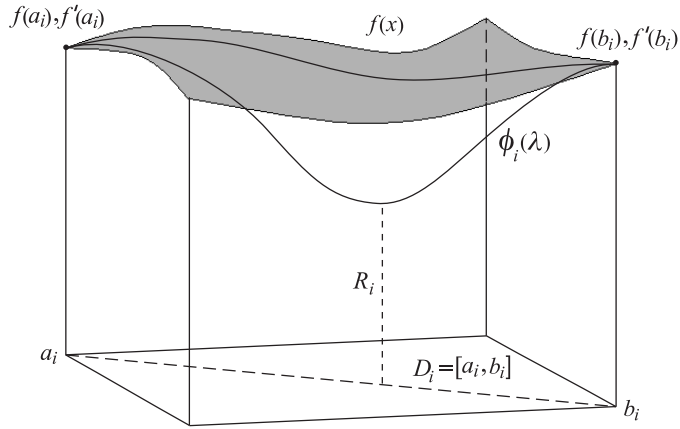}}
\caption{Obtaining the lower bound $R_i$ for the objective
function $f(x)$ with the Lipschitz gradient over $D_i = [a_i,
b_i]$ by using smooth auxiliary function~$\phi_i(\lambda)$,
$\lambda \in [a_i, b_i]$, along the main diagonal $[a_i, b_i]$ of
$D_i$.} \label{fig:DiagonalDerivativesSmooth}
\end{figure}

Regarding the partition strategies (see partitioning operator $P$
on Step~4 in Figure~\ref{fig:DBA}), the \emph{diagonal} partition
strategies (see the references in \cite{Pinter(1996),
Sergeyev&Kvasov(2006), Sergeyev&Kvasov(2008),
Sergeyev&Kvasov(2011)}) are taken into account in this paper,
although other types of subdivisions are worth to be considered
(see, e.\,g., simplicial partitions in
\cite{Clausen&Zilinskas(2002), Paulavicius:et:al.(2014),
Paulavicius&Zilinskas(2014a)}). In the diagonal approach, the
hyperinterval $D$ from (\ref{LGOP_D}) is subdivided into a set of
smaller hyperintervals, the objective function (and its gradient)
is (are) evaluated only at two vertices corresponding to the main
diagonal of the generated hyperintervals (see, e.\,g., vertices
$a_i$ and $b_i$ of a hyperinterval $D_i$ in
Figure~\ref{fig:DiagonalDerivativesSmooth}), and the results of
these trials are used to choose a hyperinterval for subsequent
subdivisions. This approach has a number of interesting
theoretical properties and has proved to be efficient in solving
practical problems.

First, efficient one-dimensional global optimization algorithms
can be easily extended to the multidimensional case by means of
the diagonal approach (see, e.\,g., \cite{Sergeyev&Kvasov(2006),
Sergeyev&Kvasov(2008), Sergeyev&Kvasov(2011)}). For instance, in
order to obtain the characteristic $R_i$ of a multidimensional
hyperinterval $D_i$, some univariate characteristics can be used
as prototypes. After an appropriate transformation they can be
considered along the one-dimensional segment being the main
diagonal $[a_i, b_i]$ of the hyperinterval $D_i$ (see smooth
auxiliary function $\phi_i(\lambda)$ in
Figure~\ref{fig:DiagonalDerivativesSmooth}).

Second, the diagonal approach is computationally close to one of
the simplest strategies---center-sampling technique (see,
e.\,g.,~\cite{Conn:et:al.(2009), Evtushenko(1985),
Evtushenko&Posypkin(2011), Gablonsky&Kelley(2001),
Jones:et:al.(1993)})---but at the same time, the function trials
are executed at two points of each hyperinterval, providing so
more information about the function over the subregion under
investigation than center-sampling algorithms that evaluate~$f(x)$
(and~$f'(x)$) only at one central point.

An important issue in the framework of diagonal algorithms is the
way the hyperintervals are partitioned during the search. Various
exploration techniques based on different diagonal partition
strategies are studied, e.\,g., in~\cite{Sergeyev(2000),
Sergeyev&Kvasov(2008), Sergeyev&Kvasov(2011)}. It is shown that
partition strategies traditionally used in diagonal methods can
perform many redundant trials and, hence, do not fulfil the
requirements of computational efficiency (this redundancy
significantly slows down the search in the case of expensive
functions).

\begin{figure}[ht]
\begin{center}
\includegraphics[width=0.8\textwidth,keepaspectratio]{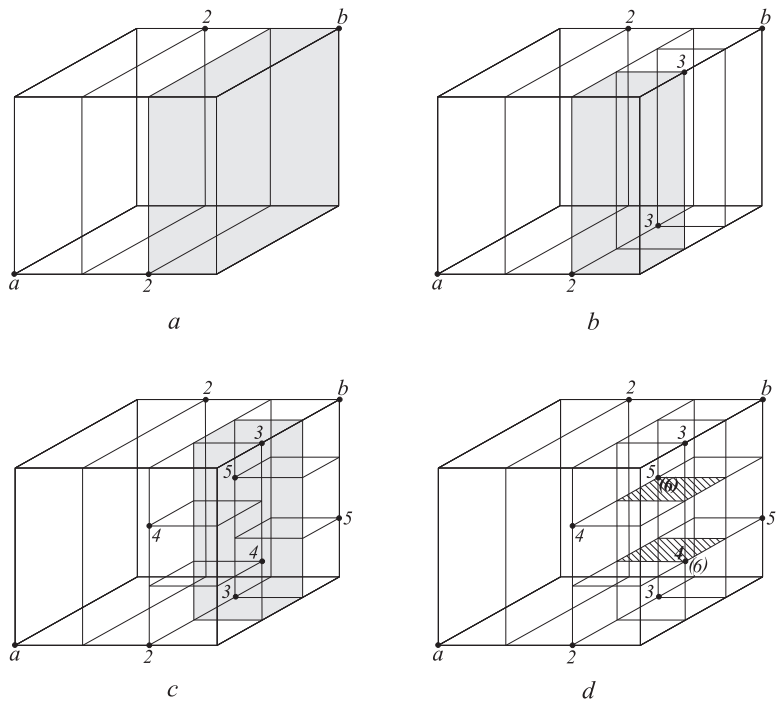}
\caption{Example of partitioning a three-dimensional hyperinterval
by the efficient diagonal partition strategy.}
\label{fig:NewStrategy_3D}
\end{center}
\end{figure}

A diagonal partition strategy introduced in~\cite{Sergeyev(2000),
Sergeyev&Kvasov(2008)} allows one to avoid the computational
redundancy of diagonal methods. In contrast to traditional
diagonal schemes, the proposed strategy (called
\emph{non-redundant}, or \emph{efficient} diagonal partition
strategy) produces regular trial meshes in such a way that one
vertex where~$f(x)$ and $f'(x)$ are evaluated can belong to
several hyperintervals (up to $2^N$, $N$~is the problem dimension
from~(\ref{LGOP_D}); see Figure~\ref{fig:NewStrategy_3D} that
shows how partitions of gray-colored hyperintervals are executed
in the course of iterations 1--6, where trial points are
represented by black dots). So, the time-consuming operation of
evaluating both $f(x)$ and $f'(x)$ is replaced by a significantly
faster procedure of reading (up to $2^N$ times) the previously
obtained functions values, saved in a special database (see,
e.\,g., \cite{Sergeyev&Kvasov(2008)}). In Figure
\ref{fig:NewStrategy_3D}, during the 6-th iteration the partition
has been performed \emph{for free} since the required values of
$f(x)$ and $f'(x)$ have been already obtained at the 4-th and 5-th
iterations (trial points of the 6-th iteration are in brackets,
see Figure \ref{fig:NewStrategy_3D}d). Hence, the non-redundant
partition strategy speeds up the search and also leads to saving
computer memory. It is important that this feature becomes more
pronounced when the problem dimension $N$ increases (see, e.\,g.,
\cite{Kvasov&Sergeyev(2003), Sergeyev&Kvasov(2006),
Sergeyev&Kvasov(2008)}). A formal procedure of subdivisions by
means of the efficient diagonal partition strategy can be viewed
in the next Section (see formulae
(\ref{u_LTDiagSmooth})--(\ref{Dm+2_LTDiagSmooth})).

A novel scheme for developing fast Lipschitz global optimization
methods can be, thus, considered. It can rely on the non-redundant
diagonal partition strategy allowing an efficient extension of
well-studied one-dimensional Lipschitz global optimization
algorithms to the multidimensional case. Interesting
multidimensional diagonal methods, based on different ways for
achieving the Lipschitz information and constructed in the
framework of the efficient diagonal approach, have been proposed
and their convergence properties have been studied, e.\,g.,
in~\cite{Kvasov&Sergeyev(2003), Sergeyev&Kvasov(2006),
Sergeyev&Kvasov(2008)}. In the following Section, this scheme is
used (see~\cite{Sergeyev&Kvasov(2008)}) to develop a new global
optimization algorithm for solving problem
(\ref{LGOP_f})--(\ref{LGOP_D}). This algorithm uses smooth
diagonal auxiliary functions along the main diagonals of
hyperintervals and shows how the unknown constant $K$ from
(\ref{LGOP_K}) can be adaptively estimated during the search.

\section{New algorithm {\sc SmoothD}}

In this Section, a new algorithm {\sc SmoothD} for solving problem
(\ref{LGOP_f})--(\ref{LGOP_D}) is described. It generalizes the
one-dimensional method adaptively estimating the constant $K$
(see~\cite{Sergeyev(1998a)}) to the multidimensional case by using
the diagonal approach and the efficient partition strategy
from~\cite{Sergeyev(2000), Sergeyev&Kvasov(2008)}. First, the new
method is presented and its computational scheme is given, then
its convergence properties are analyzed.

To describe the algorithm {\sc SmoothD} formally, we need the
following designations:

$k \geq 1$ -- the current iteration number of the algorithm;

$p(k)$ -- the total number of function trials executed during the
previous $k-1$ iterations of the method;

$\{x^{p(k)}\}$ -- the sequence of trial points generated by the
algorithm in $k$ iterations;

$\{z^{p(k)}\} = \{f(x^{p(k)})\}$ -- the corresponding sequence of
the objective function values;

\rule{0mm}{7mm}$\{ \frac{\partial z^{p(k)}}{\partial x(j)}
\}=\{\frac{\partial f(x^{p(k)})}{\partial x(j)}\}$, $1 \leq j \leq
N$, -- $N$ sequences of the corresponding partial derivatives
$f'_j(x)$, $1 \leq j \leq N$, evaluated at $x^{p(k)}$;

$M = M(k)$ -- the total number of hyperintervals of the current
partition of the search domain at the beginning of the $k$-th
iteration;

$\{D^k\}  = \{D_i^k\} = \{D_i\}$, $1 \leq i \leq M$, -- the
current partition of the initial hyperinterval $D$
into hyperintervals $D_i = [a_i,b_i]$.

For each hyperinterval $D_i$ the coordinates of the vertices $a_i$
and $b_i$ of its main diagonal $[a_i, b_i]$, the corresponding
values of $f(a_i)$ and $f(b_i)$, and the corresponding values of
partial derivatives $f'_j(a_i)$ and $f'_j(b_i)$, $1 \leq j \leq
N$, are stored.

Directional derivatives along main diagonals are then evaluated as follows:

 \be
  f'(a_i) = \left(\sum_{j=1}^{N} f'_j(a_i)(b_i(j)-a_i(j))\right)/\Delta_i,
  \label{dz_ai}
 \ee
 \be
  f'(b_i) = \left(\sum_{j=1}^{N} f'_j(b_i)(b_i(j)-a_i(j))\right)/\Delta_i,
  \label{dz_bi}
 \ee
where
 \be
  \Delta_i = \| a_i - b_i \| = \sqrt{\sum_{j=1}^{N} (a_i(j) -
  b_i(j))^2} \label{kniga5_Delta_i}
 \ee
is the length of the main diagonal of hyperinterval $D_i$, $1 \leq i \leq
M(k)$.

{\bf New Algorithm {\sc SmoothD}.} Set the initial value of the
iteration counter $k:=1$. Evaluate the objective function $f(x)$
and its partial derivatives $f'_j(x)$, $1 \leq j \leq N$, at the
vertices~$a$ and~$b$ of the search domain~$D$: $x^1 = a$, $x^2 =
b$. Set both the current number of trials $p(1):=2$ and the
current number of hyperintervals $M(1):=1$. The initial partition
of the search domain is~$\{D^1\} := \{D\}$.

Suppose that $k \geq 1$ iteration of the algorithm have been already executed. Iteration $k+1$ consists of the following Steps.

\begin{description}

\item {\bf Step 1} (\textit{Estimation of the Lipschitz constant for the gradient}).
Obtain the current estimate $m$ of the Lipschitz constant $K$ from
(\ref{LGOP_K}) for the objective function gradient as
 \be
   m  = r \max \{ \xi, \hat{m} \}, \label{K_LTDiagSmooth}
 \ee
where $r>1$ is the reliability parameter of the method, $\xi$ is a
small positive value (it ensures the correct algorithm execution
when the value $\hat{m}$ is too small; as a rule, $\xi \leq
10^{-6}$), and $\hat{m}$ is found as
 $$
   \hat{m} = \max_{1 \leq i \leq M(k)} \{ w_i\},
 $$
with
 $$
  w_i=\frac{|2(f(a_i)-f(b_i))+(f'(a_i)+f'(b_i))\Delta_i| +
  d_i}{\Delta_i^2}
 $$
and
 \begin{multline*}
   d_i=\{[2(f(a_i)-f(b_i))+(f'(a_i)+f'(b_i))\Delta_i]^2 \,+ \\
   +\, (f'(b_i)-f'(a_i))^2\Delta_i^2\}^{\frac{1}{2}}.
 \end{multline*}

\item {\bf Step 2} (\textit{Characteristics calculations}). For each hyperinterval \mbox{$D_i \in \{D^k\}$}, $1 \leq i \leq
M(k)$, calculate its characteristic~$R_i$, being the minimum of
the auxiliary function $\phi_i(\lambda)$ on the diagonal $[a_i,
b_i]$ (see Fig. \ref{fig:DiagonalDerivativesSmooth}), as follows.

 {\bf Step 2.1.} If inequality
 $$
  \pi'_{i}(y_{i}) \pi'_{i}(y'_{i}) < 0,
 $$
is verified, where $\pi'_i (y_{i})$, $\pi'_i (y'_{i})$ are found
by formulae
 $$
  \pi'_{i}(y_{i})  = m y_{i} + B_i, \hspace{5mm}
  \pi'_{i}(y'_{i}) = m y'_{i} + B_i,
 $$
with
 \begin{multline*}
 {\displaystyle y_i=\frac{\Delta_i}{4}+\frac{f'(b_i)-f'(a_i)}{4m} \,+ } \\
 {\displaystyle +\, \frac{f(a_i)-f(b_i)+f'(b_i)\Delta_i + 0.5m\Delta_i^2}{m\Delta_i+f'(b_i)-f'(a_i)},}
 \end{multline*}
 \vspace*{-3mm}
 \begin{multline*}
 {\displaystyle y_i' = -\frac{\Delta_i}{4}-\frac{f'(b_i)-f'(a_i)}{4m} \,+}\\
 {\displaystyle +\,\frac{f(a_i)-f(b_i)+f'(b_i)\Delta_i +
0.5m\Delta_i^2}{m\Delta_i+f'(b_i)-f'(a_i)},}
 \end{multline*}

 $$
  B_i = f'(b_i) - 2m y_i + m\Delta_i,
 $$
then go to~Step~2.2. Otherwise, go to~Step~2.3.

 {\bf Step 2.2.} Calculate characteristic $R_i$ as
 \be \label{Ri_1_LTDiagSmooth}
   R_i = \min \{f(a_i), \phi_i(\hat{x_i}), f(b_i) \},
 \ee
where
 $$
   \phi_i(\hat{x_i}) = f(b_i) - f'(b_i)\Delta_i - 0.5m\Delta_i^2 +
   m y_i^2 - 0.5m\hat{x_i}^2
 $$
and
 $$
  \hat{x_i} = 2y_i - m^{-1}f'(b_i) - \Delta_i.
 $$
Go to~Step~3.

{\bf Step 2.3.} Calculate characteristic $R_i$
as
 \be \label{Ri_2_LTDiagSmooth}
   R_i = \min \{f(a_i), f(b_i) \}.
 \ee

\item {\bf Step 3} (\textit{Hyperinterval selection}). Select a hyperinterval~\mbox{$D_t \in
\{D^k\}$} for the further partitioning such that condition
 \be \label{minR_LTDiagSmooth}
  t = \arg\min_{1 \leq i \leq M(k)} {R_i}
 \ee
is verified.
If there are several hyperintervals satisfying condition (\ref{minR_LTDiagSmooth}), then choose among them the hyperinterval with the minimal index $t$.

\item {\bf Step 4} (\textit{Stopping criterion}). If
  \be
   \| a_t - b_t \| \leq \varepsilon \| a - b
   \|,  \label{Stop_LTDiagSmooth}
  \ee
where $\varepsilon >0$ is a given accuracy, then stop the algorithm.

Take as an estimate of the global minimum~$f^*$
value $z_k^*$
 \[
  z^*_k=\min \{z \mid z \in z^{p(k)}\},
 \]
attained at point
 \[
  x^*_k={\rm arg} \min_{x^i \in x^{p(k)}} f(x^i).
 \]
Otherwise, go to Step~5.

\item {\bf Step 5} (\textit{Generation of points~$u$ and~$v$}).
Calculate coordiantes of points~$u$ and~$v$ for partitioning the hyperinterval $D_t=[a^k_t, b^k_t]$:
 \begin{multline} \label{u_LTDiagSmooth}
  u=(a(1),\ldots,a(i-1),a(i) \,+ \\
  +\, \frac{2}{3} (b(i)-a(i)),a(i+1),\ldots,a(N)),
 \end{multline}
 \vspace*{-5mm}
 \begin{multline}\label{v_LTDiagSmooth}
  v=(b(1),\ldots,b({i-1}),b(i) +\, \\
  +\frac{2}{3}(a(i)-b(i)),b(i+1),\ldots,b(N)),
 \end{multline}
where $a(j)=a^k_t(j),\ b(j)=b^k_t(j),\ 1 \leq j \leq N$, and~$i$ is found as
 \be \label{maxD_LTDiagSmooth}
   i = \min\, \arg \max_{1 \leq j \leq N}\{ | b(j) - a(j) | \}.
 \ee

Check whether the function $f(x)$ and its partial derivatives
$f_j'(x)$, $1 \leq j \leq N,$ have been already evaluated at the
points $u$ and~$v$.

{\bf Step 5.1.} If $f(x)$ and~$f_j'(x)$, $1 \leq j \leq N$, have been previously evaluated at both the points $u$ and~$v$, then set $p(k+1):=p(k)$ and go to Step~6.

{\bf Step 5.2.} If $f(x)$ and~$f_j'(x)$, $1 \leq j \leq N$, have been previously evaluated at only one point, e.\,g., at point $u$ (similar  operation is performed in the case when the functions evaluations have been performed only at point $v$), then perform a new trial of $f(x)$ and~$f_j'(x)$, $1 \leq j \leq N$, at point
 \[
  x^{p(k)+1} = v\,,
 \]
set $p(k+1) := p(k)+1$ and go to Step~6.

{\bf Step 5.3.} If neither $f(x)$ nor~$f_j'(x)$, $1 \leq j \leq N$, have been previously evaluated at points~$u,\,v$, then evaluate $f(x)$ and~$f_j'(x)$, $1 \leq j \leq N$, at points
 \[
   x^{p(k)+1} = u, \ \ x^{p(k)+2} = v,
 \]
set $p(k+1) := p(k)+2$ and go to Step~6.

\item {\bf Step 6} (\textit{Efficient diagonal partition}).
Obtain a new partition $\{D^{k+1}\}$ of the search domain $D$:
 \be  \label{artitionD_LTDiagSmooth}
   {\displaystyle \{D^{k+1}\}= \{D^{k}\} \setminus \{D_t\}\cup
   \{D^{k+1}_{M(k)+1}\} \cup \{D^{k+1}_{M(k)+2}\},}
 \ee
where vertexes of the main diagonals of new
hyperintervals~$D^{k+1}_t$, $D^{k+1}_{M(k)+1}$
and~$D^{k+1}_{M(k)+2}$ are given by formulae
 \be
   a^{k+1}_t = a_{M(k)+2} = u, \hspace{2mm}
   b^{k+1}_t = b_{M(k)+1} = v. \label{Dt_LTDiagSmooth}
 \ee
 \be
   a_{M(k)+1} = a^k_t, \hspace{2mm}
   b_{M(k)+1} = v, \label{Dm+1_LTDiagSmooth}
 \ee
 \be
   a_{M(k)+2} = u, \hspace{2mm}
   b_{M(k)+2} = b^k_t, \label{Dm+2_LTDiagSmooth}
 \ee
respectively.

Set $M(k+1):= M(k)+2$, increase the iteration counter
$k:=k+1$ and go to Step~1.

\end{description}

The introduced algorithm {\sc SmoothD} belongs to the class of
diagonally extended geometric algorithms and, thus, to a more
general class of `Divide-the-Best' algorithms (see
\cite{Sergeyev(1998b)}). Its reliability parameter $r$ from
(\ref{K_LTDiagSmooth}) controls the estimate $m$ of the Lipschitz
constant $K$ of $f'(x)$ and affects the algorithm's convergence,
as it will be shown in the next Section. Namely, by increasing $r$
the reliability of the method also increases. As this parameter
decreases, the search rate increases, but the probability of
convergence to a point other than the global minimizer of $f(x)$
grows as well.

The next Theorem establishes convergence properties of the proposed algorithm.

 \begin{theorem}
For any function $f(x)$ with the gradient satisfying the Lipschitz
condition~$(\ref{LGOP_K})$ with the constant $K,\,\,0<K<\infty,$
there exists a value~$r^*$ of the reliability parameter $r$ of the
{\sc SmoothD} algorithm such that for any $r \geq r^*$ the
infinite ($\varepsilon = 0$ in (\ref{Stop_LTDiagSmooth})) sequence
of trial points $\{x^{p(k)}\},$ generated by the {\sc SmoothD}
algorithm during minimization of $f(x),$ will converge only to the
global minimizers of $f(x)$.
 \label{DGG_BS_r*}
 \end{theorem}

{\bf Proof.} This result can be obtained as a particular case of
the general convergence study of `Divide-the-Best' algorithms from
\cite{Sergeyev(1998b)} and its proof is so omitted. $\hfill \Box $

It should be noticed in this context that if a value of $r$
smaller than~$r^*$ is used, the algorithm can converge (see the
general analysis executed for `Divide-the-Best' methods in
\cite{Sergeyev(1998b)}) to a local minimizer of $f(x)$ or to the
boundary of a subregion of $D$ corresponding to the best
characteristics~(\ref{Ri_1_LTDiagSmooth})--(\ref{minR_LTDiagSmooth}).
This situation indicates the necessity to increase the value of
the reliability parameter, i.\,e., it is a practical hint for the
choice of $r$.

\section{Numerical experiments}

In this Section, results of numerical experiments aimed at
evaluating the {\sc SmoothD} method performance and comparing it
with some known global optimization techniques are given.
Particularly, two popular global optimization algorithms
(downloadable from
\url{http://www4.ncsu.edu/~ctk/SOFTWARE/DIRECTv204.tar.gz}) have
been taken for numerical comparisons: the DIRECT method
from~\cite{Jones:et:al.(1993)} and its locally-biased version
DIRECT{\it l} from \cite{Gablonsky&Kelley(2001)}. Both of them
partition the search region $D$ into small hyperintervals and are
widely used for solving applied black-box global optimization
problems, being therefore important competitors for the {\sc
SmoothD} method. They also belong to the class of
`Divide-the-Best' algorithms but do not use the information about
the objective function gradient during their work.

All experiments were performed on a 3.40 GHz Intel(R) Core(TM)
i7-2600 CPU PC with 12 Gb memory using 64-bit Windows system.

This Section is structured as follows. First, methodology of
executing numerical experimentation is described presenting test
functions and comparison criteria used. Then, numerical issues on
the reliability parameter of the new method are discussed and
sensitivity analysis of this parameter is performed. Finally,
overall numerical results are summarized and commented~on.

\subsection{Classes of test functions and comparison criteria}

Eight classes of differentiable test problems generated by the
GKLS-generator from~\cite{Gaviano:et:al.(2003)} were used to
perform numerical experiments. This generator constructs three
types (non-differentiable, continuously differentiable, and twice
continuously differentiable) of classes of multidimensional and
multiextremal test functions with known local and global minima.
The generation procedure consists of defining a convex quadratic
function systematically distorted by polynomials.

Each test class provided by the generator consists of 100
functions and is defined by the following parameters: (i) problem
dimension $N$, (ii) number of local minima $m$, (iii) global
minimum value $f^*$, (iv) radius of the attraction region of the
global minimizer $\rho$, (v) distance from the global minimizer to
the quadratic function vertex $d$. The other necessary parameters
are chosen randomly by the generator for each test function of the
class. A special notebook with a complete description of all
functions is supplied to the user. The GKLS-generator always
produces the same test classes for a given set of the user-defined
parameters, allowing one to perform repeatable numerical
experiments. By changing the user-defined parameters, classes with
different properties can be created. The generator is available on
the ACM Collected Algorithms (CALGO) database (the CALGO is part
of a family of publications produced by the Association for
Computing Machinery) and it is also downloadable for free from
\url{http://wwwinfo.dimes.unical.it/~yaro/GKLS.html}.

In order to obtain comparable results, the same eight GKLS test
classes of continuously differentiable functions of dimensions $N
= 2, 3, 4$, and 5 defined by the same five parameters as
in~\cite{Sergeyev&Kvasov(2006)} were used (see Table
\ref{tab:description}). The global minimum value $f^*$ was equal
to $-1.0$ and the number of local minima $m$ was equal to 10 for
all classes (these values are default settings of the
GKLS-generator). Two test classes were considered for each
particular dimension $N$: the `simple' class and the `hard' one.
The difficulty of a class was increased either by decreasing the
radius~$\rho$ of the attraction region of the global minimizer
$x^*$ (as for two- and five-dimensional classes) or by increasing
the distance $d$ from $x^*$ to the quadratic function vertex
(three- and four-dimensional classes).

\begin{table}
\caption{Description of the GKLS test classes used in numerical
experiments} \label{tab:description}
\begin{tabular*}{\textwidth}{@{\extracolsep{\fill}}clcccccc}
\hline\noalign{\smallskip}
Class & Difficulty & $\varepsilon$ & $N$ & $m$ & $f^*$ & $d$ & $\rho$\\
\noalign{\smallskip}\hline\noalign{\smallskip}
1. & Simple & $10^{-4}$ & 2 & 10 & $-1.0$ & 0.90 & 0.20\\
2. & Hard   & $10^{-4}$ & 2 & 10 & $-1.0$ & 0.90 & 0.10\\
3. & Simple & $10^{-6}$ & 3 & 10 & $-1.0$ & 0.66 & 0.20\\
4. & Hard       & $10^{-6}$ & 3 & 10 & $-1.0$ & 0.90 & 0.20\\
5. & Simple & $10^{-6}$ & 4 & 10 & $-1.0$ & 0.66 & 0.20\\
6. & Hard       & $10^{-6}$ & 4 & 10 & $-1.0$ & 0.90 & 0.20\\
7. & Simple & $10^{-7}$ & 5 & 10 & $-1.0$ & 0.66 & 0.30\\
8. & Hard       & $10^{-7}$ & 5 & 10 & $-1.0$ & 0.66 & 0.20\\
\noalign{\smallskip}\hline\noalign{\smallskip}
\end{tabular*}
\end{table}

In all numerical experiments, the maximal allowed number of
function trials was set equal to $p_{max}=1\,000\,000$.  A problem
was considered to be solved by a method under examination if the
method generated a trial point~$x'$ in an
$\varepsilon$-neighborhood of the global minimizer $x^* \in
D=[a,b]$, i.\,e.,
\begin{equation}
\label{eq:StoppingRule} | x'(j) -x^*(j) | \le
\sqrt[n]{\varepsilon}(a(j) - b(j)), \quad 1\le j \le N,
\end{equation}
where $0<\varepsilon < 1$ is an accuracy coefficient (its values
are given in the third column of Table~\ref{tab:description}).
Such a type of stopping criterion is acceptable only when the
global minimizer $x^*$ is known, i.\,e., in the case of test
functions. When a real black-box objective function is minimized
and global minimization methods have an internal stopping
criterion (as that of the proposed algorithm, see formula
(\ref{Stop_LTDiagSmooth})), they execute a number of iterations
(that can be high) after a `good' estimate of~$f^*$ has been
achieved in order to demonstrate the `goodness' of the solution
found (see, e.\,g., \cite{Horst&Pardalos(1995), Pinter(1996),
Strongin&Sergeyev(2000)}).

Since each evaluation of a real-life black-box objective function
is usually a time-consuming operation (see,
e.\,g.,~\cite{Horst&Tuy(1996), Paulavicius&Zilinskas(2014a),
Pinter(1996), Sergeyev&Kvasov(2008), Sergeyev:et:al.(2013),
Strongin&Sergeyev(2000), Zhigljavsky&Zilinskas(2008)}), the
maximal number of function trials executed by the methods until
the satisfaction of the stopping criteria (\ref{eq:StoppingRule})
was chosen as the main criterion of the comparison (other possible
comparison indicators as the number of generated hyperintervals,
the average number of trials, and so on can be also used, see,
e.\,g., \cite{Paulavicius:et:al.(2014), Sergeyev&Kvasov(2006),
Sergeyev&Kvasov(2008), Sergeyev:et:al.(2013)}).

\subsection{The reliability parameter and its sensitivity analysis}

In order to evaluate the influence of the reliability parameter
$r$ from (\ref{K_LTDiagSmooth}) on the convergence properties of
the {\sc SmoothD} method, the two-dimensional class 1 from Table
\ref{tab:description} was taken. A problem from this class was
considered to be solved by the new method if it generated a trial
satisfying condition~(\ref{eq:StoppingRule}) (if the method
stopped due to its internal stopping criterion
(\ref{Stop_LTDiagSmooth}) and condition (\ref{eq:StoppingRule})
was not verified, the problem was considered to be unsolved by
{\sc SmoothD} with a given value of the reliability parameter).

The following numbers were taken as indicators of the {\sc
SmoothD} method performance in this analysis subject to a fixed
value $\bar{r}$ of the reliability parameter $r$, equal to all 100
functions of the class:

$S(\bar{r})$ -- number of solved problems among 100 of the test class;

$p^*(\bar{r})$ -- maximal number of trials required for the method
to satisfy condition~(\ref{eq:StoppingRule}) for \emph{all}
$S(\bar{r})$ functions of the test class;

$p_{\rm avg}(\bar{r})$ -- average number of trials required for the method to satisfy condition~(\ref{eq:StoppingRule}) for all $S(\bar{r})$ functions of the test class.

Numerical results obtained by the proposed method on the
considered two-dimensional test class when the reliability
parameter varied from $\bar{r}=1.2$ up to $\bar{r} = 5.8$ are
reported in Table \ref{tab:r_overall}.

\begin{table}
\caption{Influence of the reliability parameter on the {\sc
SmoothD} method performance in the case of the two-dimensional
GKLS test class 1} \label{tab:r_overall}
\begin{tabular*}{\textwidth}{@{\extracolsep{\fill}}ccccccc}
\hline\noalign{\smallskip}
Indicator & $\bar{r}=1.2$ & $\bar{r}=1.8$ & $\bar{r}=2.8$ & $\bar{r}=3.8$ & $\bar{r}=4.8$ & $\bar{r}=5.8$\\
\noalign{\smallskip}\hline\noalign{\smallskip}
$p^*(\bar{r})$  &199    &272    &332    &410    &424    &451 \\
$p_{\rm avg}(\bar{r})$  &105.14 &169.63 &222.16 &293.52 &323.42 &341.60 \\
$S(\bar{r})$ &51     &81     &91     &98     &99     &100 \\
\noalign{\smallskip}\hline\noalign{\smallskip}
\end{tabular*}
\end{table}

It can be seen from Table \ref{tab:r_overall} that with increasing
the value of $r$ the reliability (the number of solved problems
$S(\bar{r})$) of the method increases. However, the numbers of
trials (both maximal and average) increase too. This is due to the
fact (well known in Lipschitz global optimization; see, e.\,g.,
\cite{Sergeyev:et:al.(2013)}) that every function optimized by the
proposed method has its own value $r^*$ from Theorem 1. Therefore,
when one executes tests with a class of 100 different functions it
becomes difficult to use specific values of $r$ for each function.

\begin{figure}[tp]
\centering
\begin{minipage}{.45\textwidth}
\centerline{\includegraphics[width=0.8\linewidth]{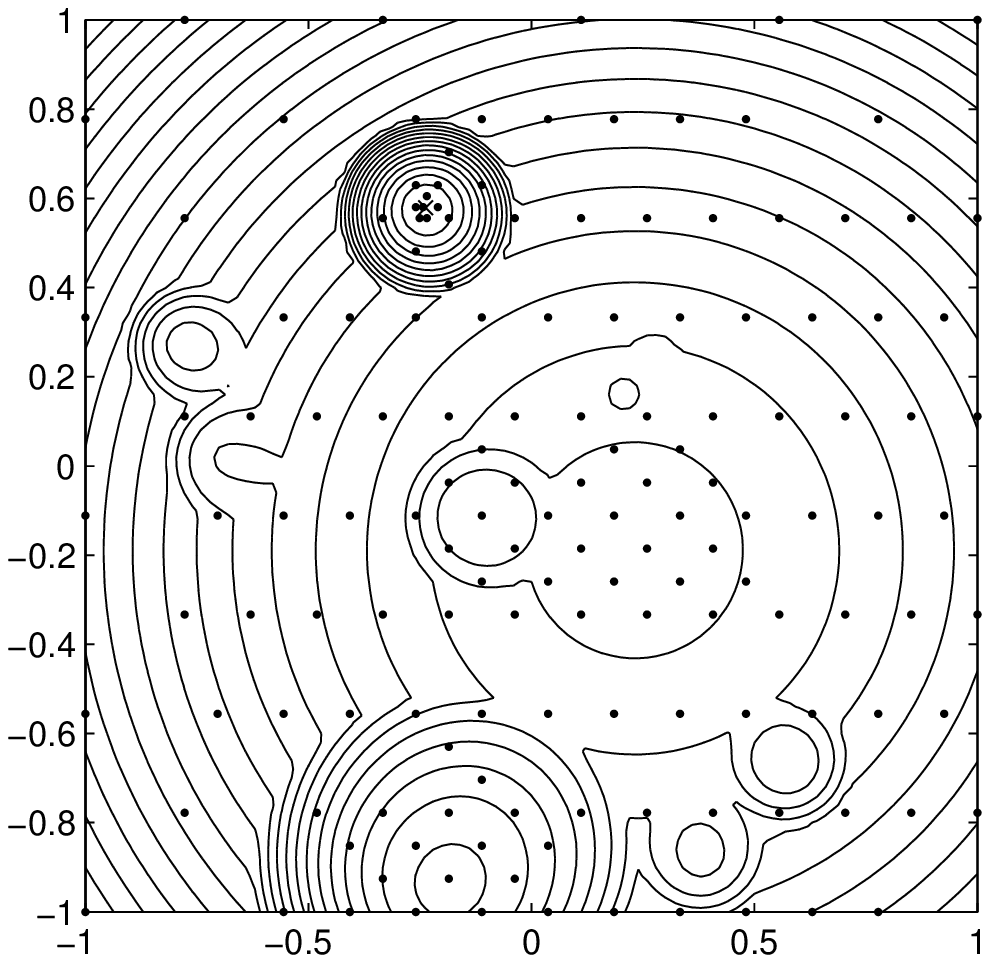}}
\end{minipage}
\begin{minipage}{.45\textwidth}
\centerline{\includegraphics[width=0.8\linewidth]{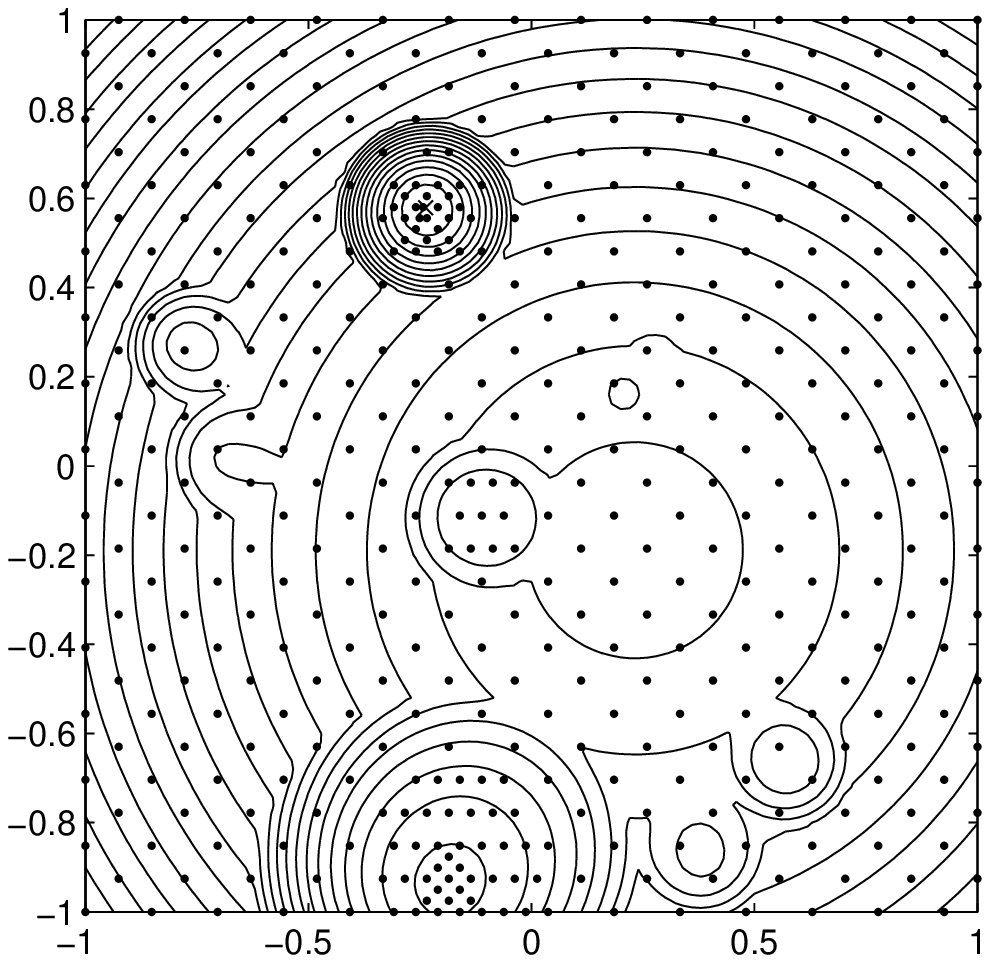}}
\end{minipage}
\caption{Trial points (black dots) generated by the {\sc SmoothD}
method when minimizing function number 58 from the GKLS test class
1, $x^* \approx (-0.2371, 0.5791)$ (\emph{left}: $\bar{r}=1.2$,
152 trials; \emph{right}: $\bar{r}=5.8$, 452 trials).} \label{f58}
\vspace*{3mm} \centering
\begin{minipage}{.45\textwidth}
\centerline{\includegraphics[width=0.8\linewidth]{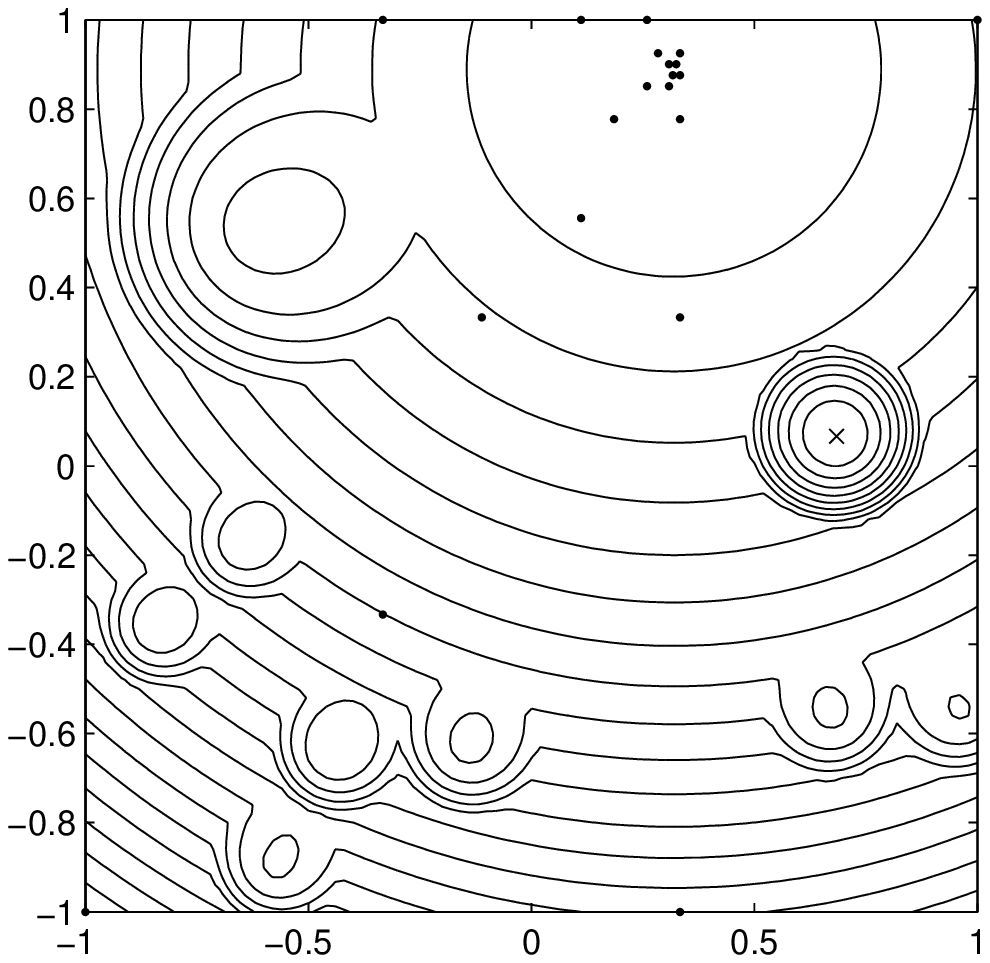}}
\end{minipage}
\begin{minipage}{.45\textwidth}
\centerline{\includegraphics[width=0.8\linewidth]{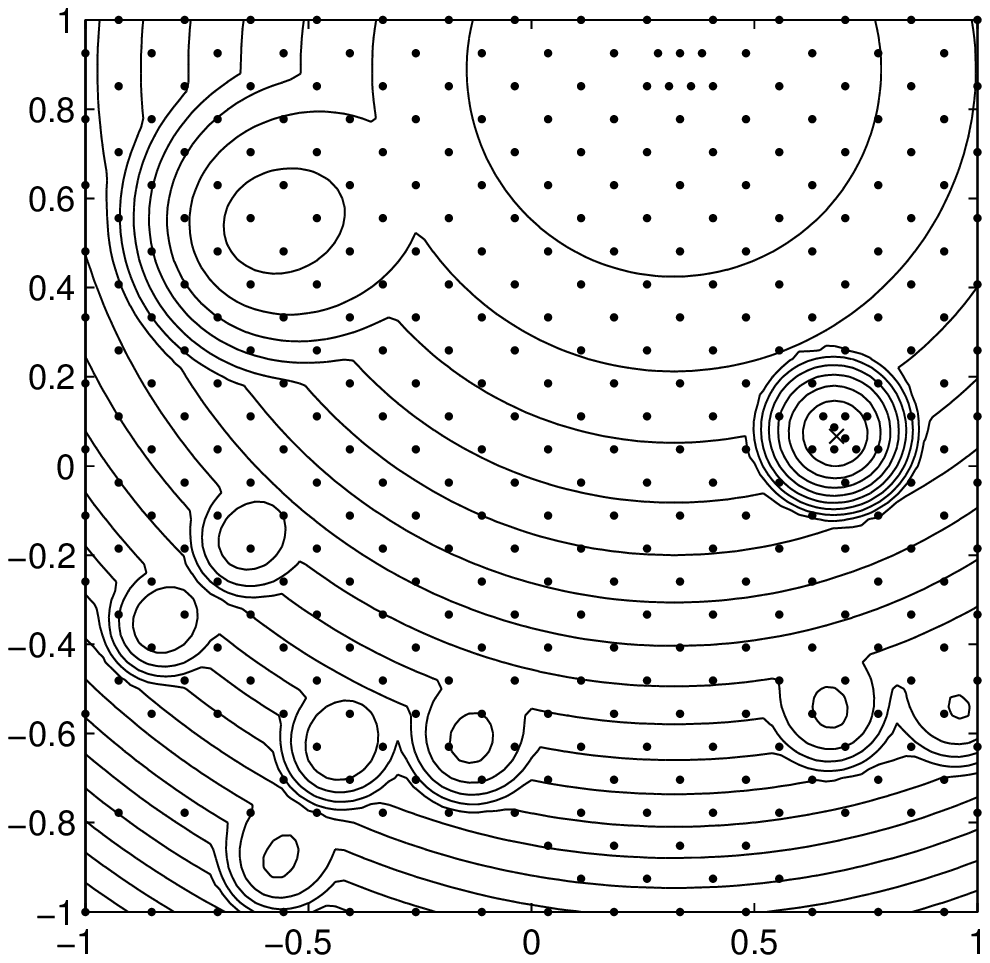}}
\end{minipage}
\caption{Trial points (black dots) generated by the {\sc SmoothD}
method when minimizing function number 54 from the GKLS test class
1 , $x^* \approx (0.6841, 0.0664)$ (\emph{left}: $\bar{r}=1.2$, 20
trials; \emph{right}: $\bar{r}=5.8$, 375 trials).} \label{f54}
\end{figure}

In fact (see an illustration given in Figures \ref{f58} and
\ref{f54}), for one particular function (see Figure \ref{f58},
left part) a relatively small value of $r$ can be already
sufficient to find its global minimizer in sense of
condition~(\ref{eq:StoppingRule}). In this case, a higher value of
the reliability parameter (see Figure \ref{f58}, right part) can
cause the method to perform an additional examination of the
search domain in the hope of capturing an eventual better minimum.
Obviously, this exploration increases the number of trials
performed by the method on this particular function (see Figure
\ref{f58}, right part) but allows the method to solve some other
problem (as that in Figure \ref{f54}) for which a greater value of
$r$ is required. An insufficient (i.e., too small) value of $r$
makes the method converge to a local minimizer, as it is shown in
Figure \ref{f54}, left part. It should be noticed in this
connection that in some practical applications the discovering of
a local minimizer can be also very useful.      Figure \ref{f54},
left part, shows that when local solutions are of an interest, the
algorithm {\sc SmoothD} can provide a fast convergence to  local
minimizers with small values of $r$.

In order to avoid an excessive exploration of the search region,
several approaches can be used. For example, an adaptive
adjustment of the reliability parameter can be performed using the
following formula
 \be \label{r_C}
   r = \bar{r} + \frac Ck,
 \ee
where $k$ is the iteration number of the method and $C$ is a
positive constant. In this way, at the initial iterations  of the
method (when $k$ is small) a greater attention is paid to the
exploration of the whole search domain $D$ while for high values
of $k$ the influence of $C$ diminishes and the reliability
parameter becomes closer to its fixed value $\bar{r}$. As it can
be seen from results in Table~\ref{tab:r_C} (where the same values
$\bar{r}$ as in Table \ref{tab:r_overall} were used but with
different values $C > 0$ in (\ref{r_C})), this adaptive technique
improves the method performance.

\begin{table}[tp]
\caption{Analysis of the adaptive reliability parameter of the
{\sc SmoothD} method in the case of the two-dimensional GKLS test
class 1} \label{tab:r_C}
\begin{tabular*}{\textwidth}{@{\extracolsep{\fill}}ccccccc}
\hline\noalign{\smallskip}
Indicator & $\bar{r}=1.2$ & $\bar{r}=1.8$ & $\bar{r}=2.8$ & $\bar{r}=3.8$ & $\bar{r}=4.8$ & $\bar{r}=5.8$\\
\noalign{\smallskip}\hline\noalign{\smallskip}
 &  &  & $C=10.0$  &  &  & \\
$p^*(\bar{r})$          &201    &277    &377    &410    &424    &453   \\
$p_{\rm avg}(\bar{r})$  &113.16 &176.06 &250.45 &294.64 &324.10 &342.01 \\
$S(\bar{r})$            &62     &86     &96     &98     &99     &100 \\
\noalign{\smallskip}\hline\noalign{\smallskip}

 &  &  & $C=20.0$  &  &  & \\
$p^*(\bar{r})$          &179    &283    &379    &411    &424    &453 \\
$p_{\rm avg}(\bar{r})$  &122.01 &181.26 &252.86 &295.51 &324.73 &342.47 \\
$S(\bar{r})$            &69     &86     &96     &99     &99     &100 \\
\noalign{\smallskip}\hline\noalign{\smallskip}

 &  &  & $C=50.0$  &  &  & \\
$p^*(\bar{r})$          &214    &293    &387    &414    &425    &454 \\
$p_{\rm avg}(\bar{r})$  &146.42 &194.98 &257.56 &299.23 &327.24 &343.83 \\
$S(\bar{r})$            &83     &92     &100    &100    &100    &100 \\
\noalign{\smallskip}\hline\noalign{\smallskip}

 &  &  & $C=100.0$  &  &  & \\
$p^*(\bar{r})$          &251    &314    &369    &416    &428    &456 \\
$p_{\rm avg}(\bar{r})$  &170.40 &212.98 &247.72 &305.07 &330.47 &345.85 \\
$S(\bar{r})$            &92     &97     &100    &100    &100    &100 \\
\noalign{\smallskip}\hline\noalign{\smallskip}

\end{tabular*}

\end{table}

As a practical recommendation for solving real-life global
optimization problems, the following procedure can be adopted. An
initial number of trials that are to be performed is fixed and the
method is started with a small value of $r$. Once the method
stopped, it is restarted on the same pool of trial points (without
re-evaluating $f(x)$ and $f'(x)$) but with a greater value of $r$.
If after repeating these restarts several times the method
converges to the same point, this point is accepted as an
approximation of the global minimizer. For more information on the
reliability parameter and other parameters of this kind in
Lipschitz global optimization see, e.\,g.,
\cite{Sergeyev:et:al.(2013), Strongin&Sergeyev(2000)}.

\subsection{Operating characteristics of the compared methods}

The overall numerical results for three compared methods on the
considered GKLS test classes are summarized in this subsection.
Since the DIRECT-type methods under examination do not have any
internal stopping criteria, they stopped either when condition
(\ref{eq:StoppingRule}) was satisfied or when~$p_{max}$ trials
were performed. During their work, these methods evaluate only the
objective function $f(x)$ from (\ref{LGOP_f}) (and do not evaluate
$f'(x)$ as the {\sc SmoothD} method does), therefore the operation
of executing a trial is less expensive in the DIRECT-type
algorithms with respect to that of the new method. The balancing
parameter $\epsilon$ of the DIRECT and DIRECT\textit{l} methods
(see~\cite{Gablonsky&Kelley(2001), Jones:et:al.(1993)}) was set
equal to $\epsilon=10^{-4}$, as recommended by many authors (see,
e.g., the review of different techniques for setting this
parameter in~\cite{Paulavicius:et:al.(2014)}).

Formula (\ref{r_C}) was used in the {\sc SmoothD} method to set
the reliability parameter $r$ from (\ref{K_LTDiagSmooth}), with
the value $\xi=10^{-6}$ in (\ref{K_LTDiagSmooth}) and the values
$C$ from~(\ref{r_C}) set in relation to the problems dimension
($C$ equal to $50$, $100$, $150$, and $200$ in the case of $N$
equal to 2, 3, 4, and 5, respectively). Moreover, several values
of $\bar{r}$ can be fixed for the entire class starting from the
initial value $\bar{r}= 1.1$ (the maximal values of this
coefficient were equal to 2.80, 5.80, 3.60, 4.30, 5.80, 6.60,
4.10, and 7.80 for the GKLS classes from 1 to 8, respectively).

The algorithms performance can be conveniently visualized by the
operating characteristics (introduced in 1978
in~\cite{Grishagin(1978)}, see,
e.\,g.,~\cite{Strongin&Sergeyev(2000)} for their English-language
description). The operating characteristics are formed by the
pairs $(p, S(p))$ where $p$ is the number of trials and $S(p)$ ($0
\leq S(p) \leq S=100$) is the number of test problems (among a set
of $S$ tests) solved by a method with less than or equal to $p$
function trials. Each pair (for increasing values of~$p$)
corresponds to a point on the plane. Higher is the graph of the
operating characteristics of a method with respect to another one,
better is the method performance on the considered set of tests.

\begin{figure}[t]
\centering
\begin{minipage}{.45\textwidth}
\centering
\includegraphics[width=\linewidth]{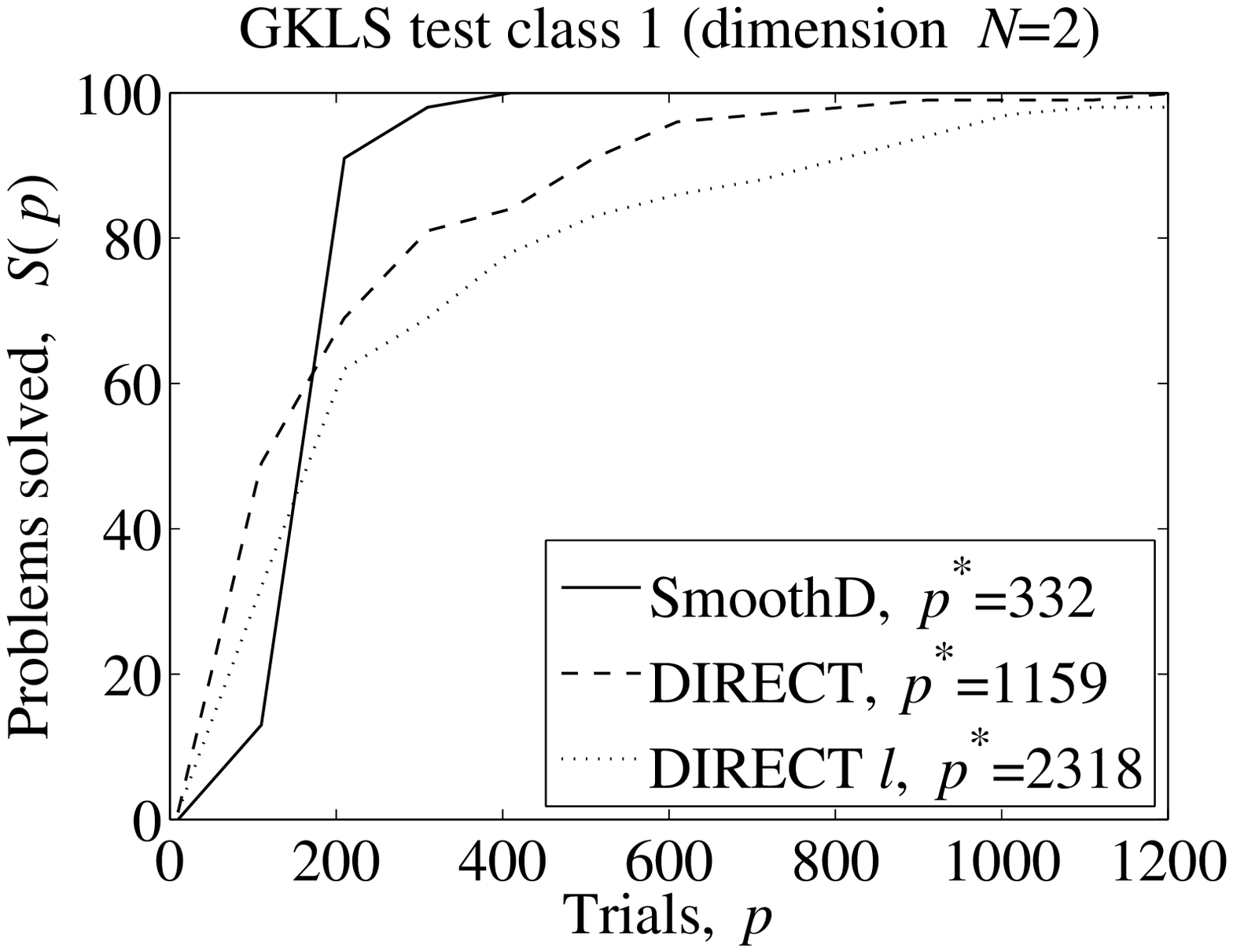}
\end{minipage}
\begin{minipage}{.45\textwidth}
\centering
\includegraphics[width=\linewidth]{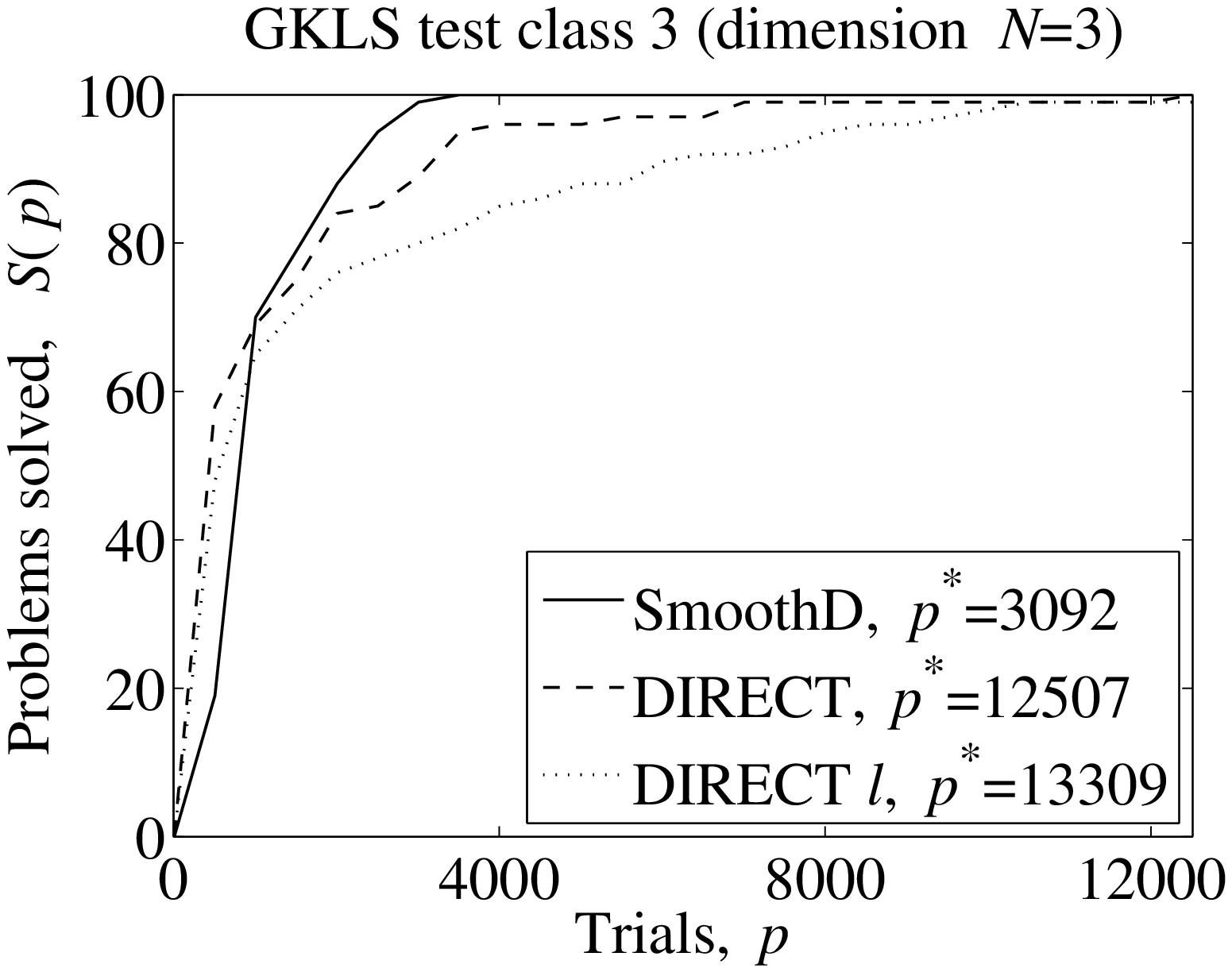}
\end{minipage}
\begin{minipage}{.45\textwidth}
\centering
\includegraphics[width=\linewidth]{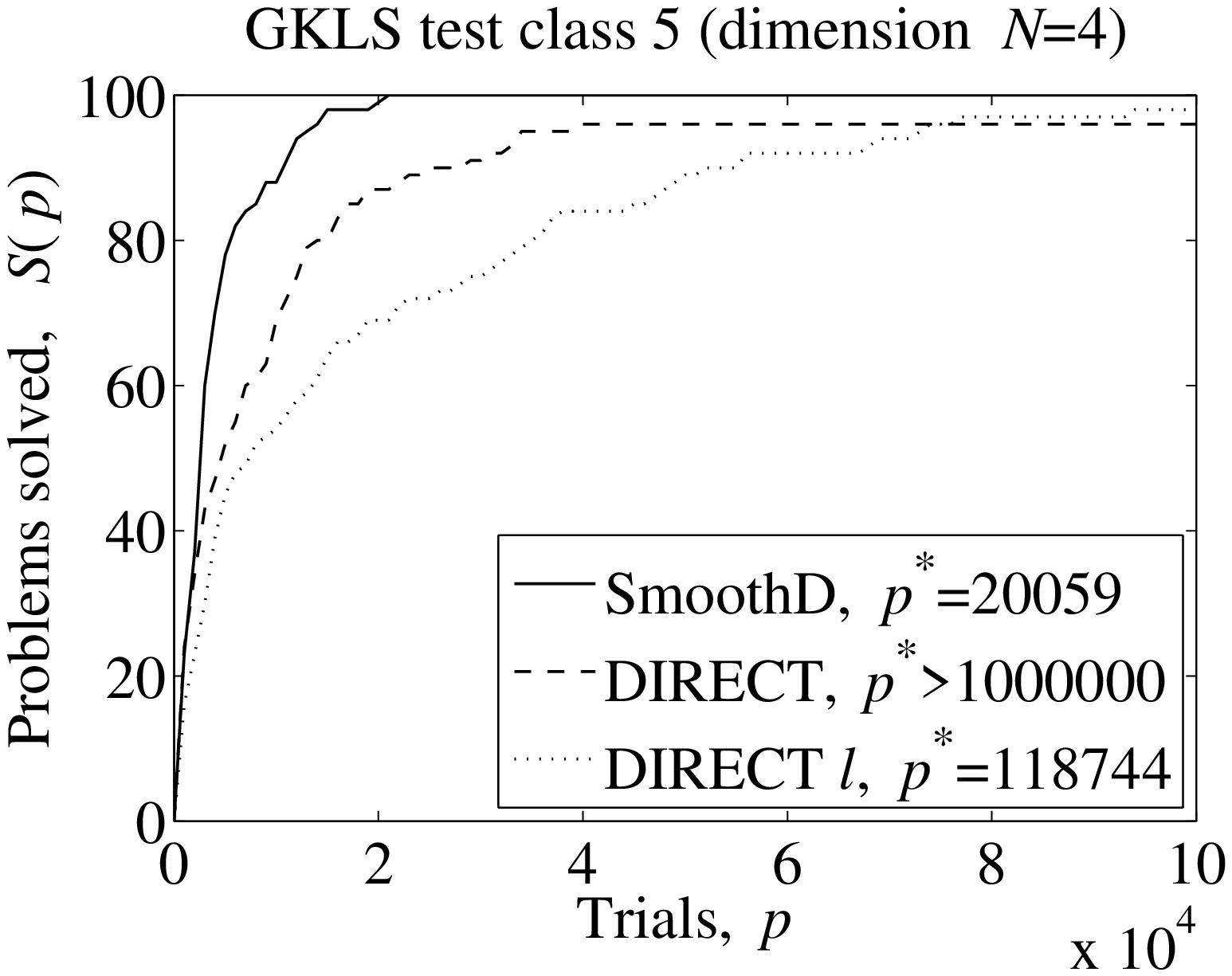}
\end{minipage}
\begin{minipage}{.45\textwidth}
\centering
\includegraphics[width=\linewidth]{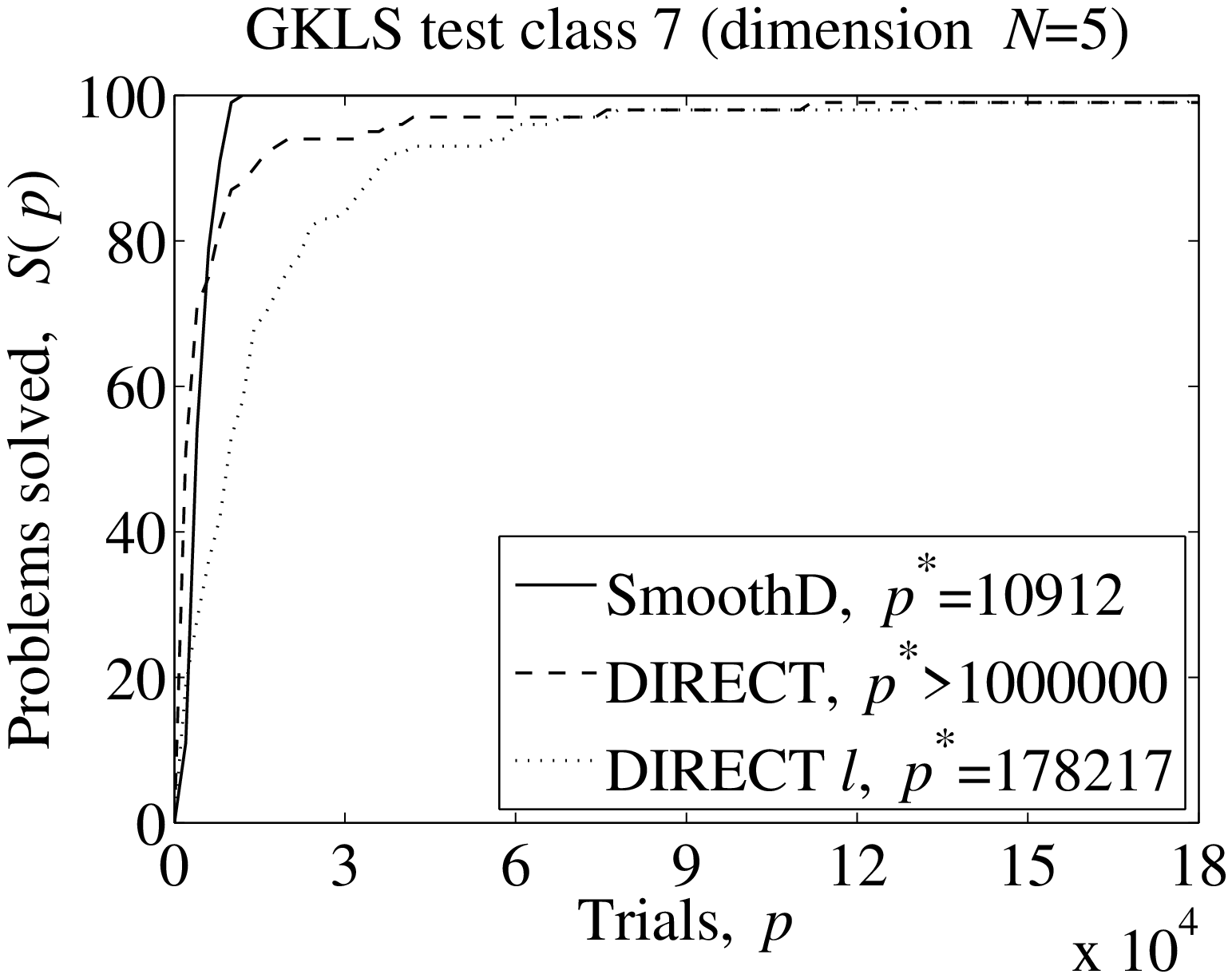}
\end{minipage}

\caption{Operating characteristics for the compared methods {\sc
SmoothD}, DIRECT, and DIRECT{\it l} on `simple' GKLS test classes
1, 3, 5, and 7 from Table \ref{tab:description}.}
\label{fig:simple}
\end{figure}

\begin{figure}[t]
\centering
\begin{minipage}{.45\textwidth}
\centering
\includegraphics[width=\linewidth]{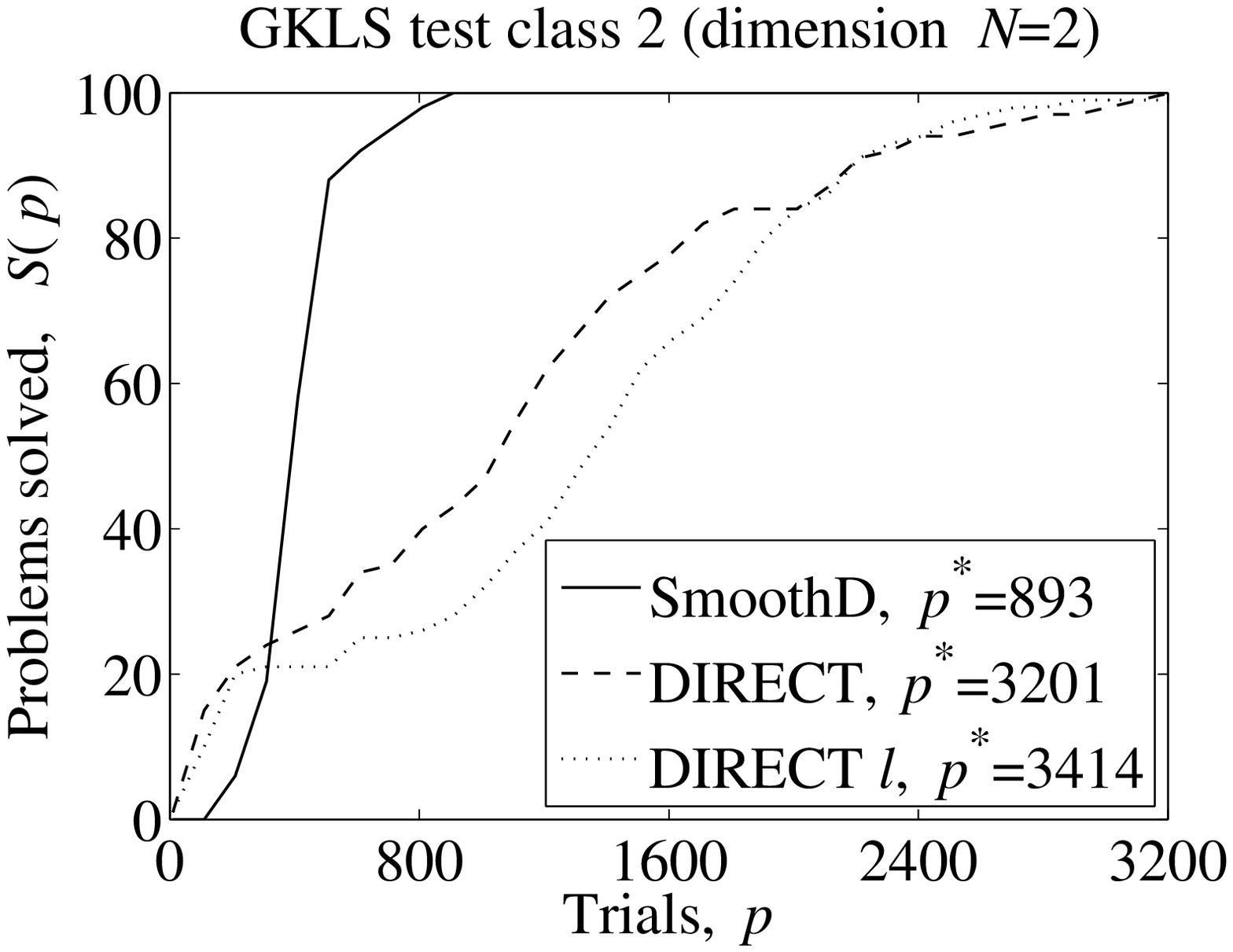}
\end{minipage}
\begin{minipage}{.45\textwidth}
\centering
\includegraphics[width=\linewidth]{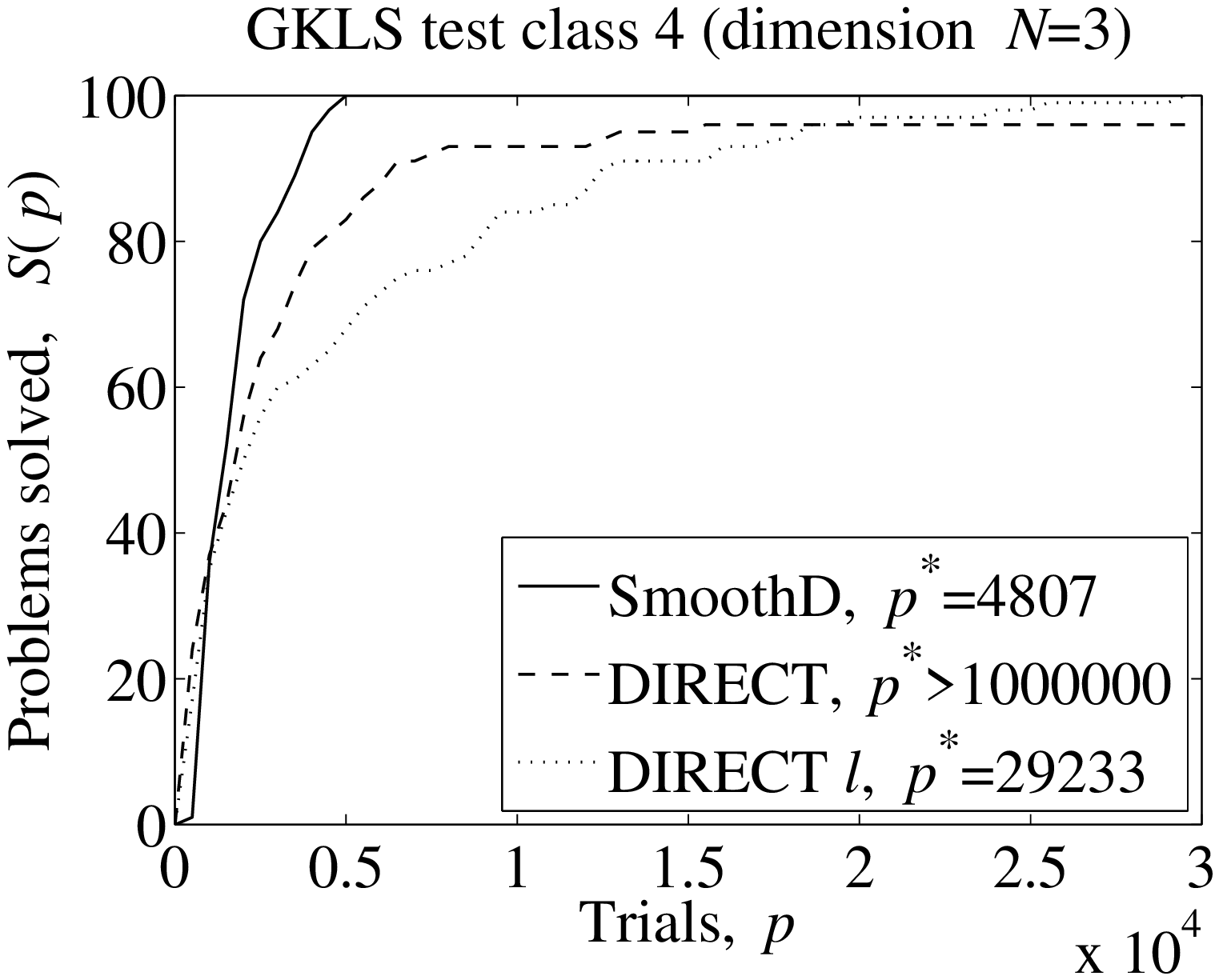}
\end{minipage}
\begin{minipage}{.45\textwidth}
\centering
\includegraphics[width=\linewidth]{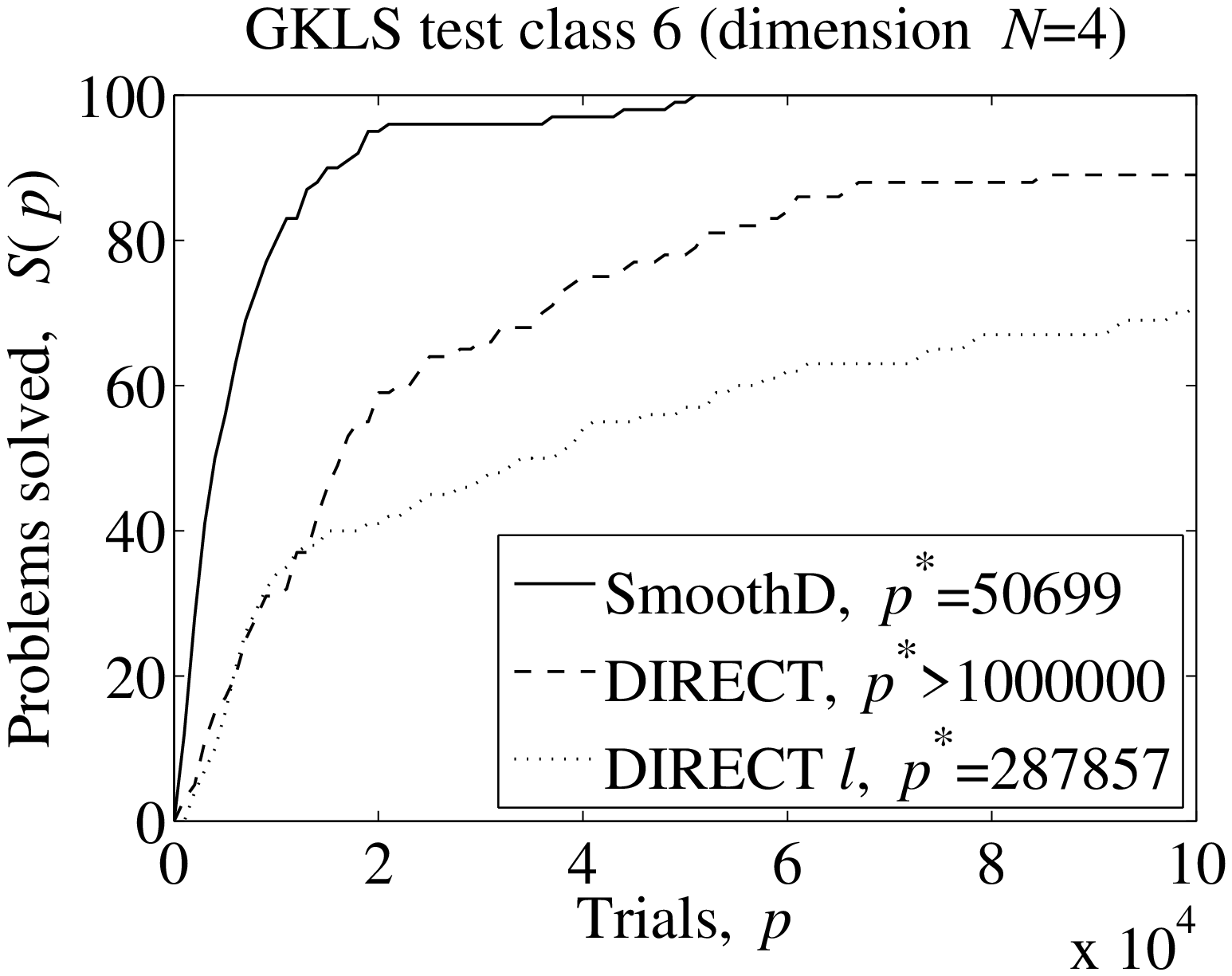}
\end{minipage}
\begin{minipage}{.45\textwidth}
\centering
\includegraphics[width=\linewidth]{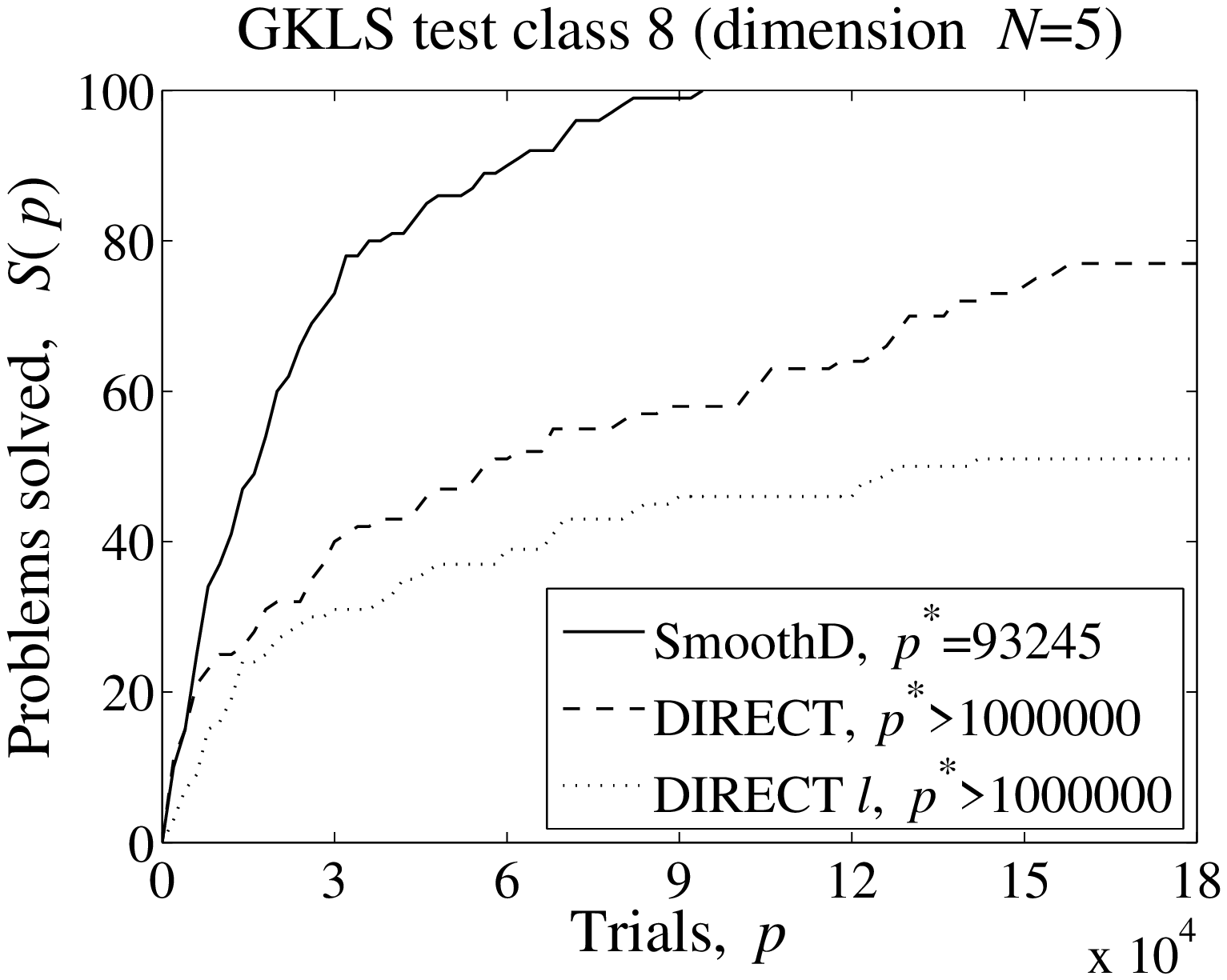}
\end{minipage}

\caption{Operating characteristics for the compared methods {\sc
SmoothD}, DIRECT, and DIRECT{\it l} on `hard' GKLS test classes 2,
4, 6, and 8 from Table \ref{tab:description}.} \label{fig:hard}
\end{figure}

In Figures \ref{fig:simple} and \ref{fig:hard}, operating
characteristics of the compared methods are shown, respectively,
on simple and hard GKLS classes from Table \ref{tab:description}.
The numbers $p^*$ of trials required for the methods to satisfy
condition~(\ref{eq:StoppingRule}) for all $100$ functions of a
test class are also indicated in these Figures (after performing
$p_{max}=1\,000\,000$ trials, the DIRECT method was unable to
solve 4, 4, 7, 1, and 16 problems of the GKLS classes 4--8,
respectively, and the DIRECT\textit{l} method was unable to solve
4 problems of the GKLS class 8).

It can be seen from the graphs reported in Figures
\ref{fig:simple} and \ref{fig:hard}, that on relatively simple
functions of each GKLS class (approximatively, half of 100
problems of the class quickly solved by every method) the {\sc
SmoothD} method behaves similarly to the DIRECT and
DIRECT\textit{l} methods, in terms of the function trials
performed (see, e.\,g., the intersection of the vertical line
corresponding to $p \approx 200$ trials and the operating
characteristics graphs on the simple two-dimensional GKLS class 1
in Figure \ref{fig:simple}). However, when more and more difficult
problems of the particular GKLS class are to be solved, the
advantage of the new method becomes more and more pronounced with
respect to the competitors (see the ascent of the operating
characteristics graphs of the {\sc SmoothD} method on all the GKLS
classes in Figures \ref{fig:simple} and \ref{fig:hard}). The most
significant advantage with respect to both the DIRECT and
DIRECT\textit{l} methods is obtained by the new method on the hard
five-dimensional GKLS class 8 (which is the most difficult among
the classes from Table \ref{tab:description}).

\section{A brief conclusion}

A new method {\sc SmoothD} has been proposed to solve global
optimization problems with the objective function having the
Lipschitz gradient. It generalizes the one-dimensional geometric
algorithm using smooth auxiliary functions and adaptive estimate
of the Lipschitz constant from~\cite{Sergeyev(1998a)} to the
multidimensional case by using the efficient diagonal scheme
from~\cite{Sergeyev(2000), Sergeyev&Kvasov(2008)}. The proposed
method belongs to the class of `Divide-the-Best' algorithms for
which strong convergence properties have been established in
\cite{Sergeyev(1998b)}. Results of numerical experiments performed
with the {\sc SmoothD} method on several hundreds of GKLS test
functions from~\cite{Gaviano:et:al.(2003)} show that already in
its basic version the proposed algorithm manifests a very nice
performance and opens new directions for developing powerful
global optimization tools. In the future, it could be interesting
to consider, in the connection to the {\sc SmoothD} method, such
approaches as the local tuning~\cite{Sergeyev(1998a),
Sergeyev&Kvasov(2008)} and local improvement
techniques~\cite{Lera&Sergeyev(2010b), Lera&Sergeyev(2013),
Sergeyev:et:al.(2013)}. Another interesting direction of research
is to develop parallel versions of the proposed method within a
general framework introduced in \cite{Gergel&Sergeyev(1999),
Sergeyev&Grishagin(1994), Sergeyev&Strongin(1989),
Strongin&Sergeyev(2000)} for constructing parallel global
optimization techniques.

\bibliographystyle{plain}
\bibliography{Deriv_Diagonal}

\end{document}